\documentclass{article}

\usepackage[latin1]{inputenc}
\usepackage{amsfonts}
\usepackage{amsmath}
\usepackage{amssymb}
\usepackage{amsthm}
\usepackage[T1]{fontenc}
\usepackage{graphicx}
\usepackage{verbatim}

\theoremstyle{definition} \newtheorem{algo}{Procedure}
\theoremstyle{definition} \newtheorem{defin}{Definition}
\theoremstyle{plain} \newtheorem{prop}{Property}[section]
\theoremstyle{plain} \newtheorem{thm}{Theorem}[section]
\theoremstyle{plain} \newtheorem{cor}{Corollary}[section]

\title{Structural properties of Stochastic Abelian Sandpile}
\author{Ayush Choure \footnote{Department of Computer Science and Engineering, IIT Bombay, Mumbai }}

\begin{document}

\maketitle

\begin{abstract}
We present some combinatorial results on the stochastic abelian sandpile model. These models are characterized by nondeterministic toppling rules. The recurrence checking for the deterministic case can be performed using the well known burning test which detects presence of forbidden sub-configurations (FSC) in strongly polynomial time. In the stochastic case, however, even for Manna's model, which is perhaps the simplest non-trivial example, no such procedure is known. In this paper, we address the decision problem of the existence of \textit{any} FSC in a general stochastic sandpile. We demonstrate a polynomial time algorithm which, given the sandpile graph and toppling rules, decides if there exists an FSC. In the event of a positive answer, it generates at least one FSC for the given sandpile. Repeated application of the algorithm can be used to find many distinct FSCs. We also demonstrate a procedure for creating larger FSCs from smaller ones and use this to create FSCs for the Manna's model. We hope that the structural analysis of stochastic sandpile we perform in this paper, will prove useful in the eventual formulation of a deterministic procedure to decide recurrence.

\end{abstract}

\section{Introduction}

The abelian sandpile model (ASM) is a type of discrete diffusion process defined on graphs with a specially designated sink vertex. The sandpile model was introduced by Bak et. al. \cite{BTW} as model of \textit{self organized criticality} in various physical systems, including the dynamics of sandpile formation. However, this model had non-commutative dynamics which were relatively less amenable to analysis. Later, Dhar \cite{DD90} introduced the abelian version of sandpile model and since then, it has been an extremely fertile area of research in the statistical physics community. It is closely associated with the \textit{chip  firing game} investigated  by  Bjorner, Lovasz and Shor \cite{auto} and Tardos \cite{chip}. The comprehensive survey article by Dhar \cite{DD06} explores its connection with diverse phenomena such as as stress distribution in earthquakes, size distribution in raindrops, path length distributions in loop-erased random walks, for instance. See Kleber \cite{gold} for a discussion on the proximity of this model with numerous disciplines such as probability theory, algorithmics, theory of computing, combinatorics, non-linear dynamics, fractals, cellular automata, to name a few.

In the deterministic abelian sandpile model, ``sand particles'' are added at the vertices of a (multi)graph. A site (vertex) is stable as long as the number of particles at the site remains less than its degree. Adding more particles would render the site unstable and is accompanied by the unstable site's passing a particle along each edge to its neighboring sites. This relaxation process is referred to as {\em{toppling}}. One of the sites known as the \textit{sink} cannot topple. To ensure that every relaxation process eventually stabilizes, one needs the condition that the sink is reachable from every other site. As the system evolves, the sandpile goes through a sequence of configurations. Those which can be revisited in any toppling sequence are called {\em{recurrent}}, the remaining ones are termed {\em{transient}}. A very basic problem here is to check if a given configuration is transient or recurrent. Dhar \cite{DD90} demonstrated a very elegant method, the \textit{burning test}, which solves the decision problem in time quadratic in graph size. The utility of the method is augmented even further as it is an essential tool in showing the existence of a bijection between a the set of spanning trees on the sandpile graph and the set of recurrent configurations \cite{DD06}.

See \cite{DD99} for an introduction to the stochastic version of the abelian sandpile model. In this model, every site has a \textit{list} of toppling rules and at instability, one of these rules is chosen arbitrarily and the toppling proceeds in accordance with that. We will define this model formally in the succeeding section. This model has a decidedly unique behavior because of the inherent non-deterministic nature. However, till now no decision procedure is known which can decide if a given configuration is transient. A transient configuration is characterized by the presence of forbidden sub-configurations (FSC), which are the minimal components that are transient by themselves. Dhar \cite{DD99} discusses some important properties of stochastic sandpile and derives bounds on the weights of recurrent configurations. A list of FSCs of Manna's model is also reported, which were discovered by empirical investigation. It is still not clear if one can determine the complete set of these configurations purely analytically.

In fact, checking for transience is equivalent to checking for presence of FSCs. However, in a general stochastic sandpile, it is not clear if an FSC will exist at all. We demonstrate the first polynomial time procedure which decides, for a given stochastic sandpile, whether there exists any FSC. Apart from solving this decision problem, our proof also sheds light on the structural aspects of sandpile graphs. We believe this result will contribute in the eventual resolution of the recurrence checking problem. We also demonstrate a simple procedure to efficiently generate arbitrarily large FSCs for the Manna's model.

\section{Abelian sandpile models}

\begin{defin}A {\em{graph}} $G$ is an ordered pair $(V(G),E(G))$ where $V(G)$ is called the set of vertices and $E(G)$ is a set of $2-$subsets of $V$, possibly with repeated elements, the set of edges.\end{defin}

This is referred to as a multi-graph in literature but we will use graph for brevity. The \textit{degree} of a vertex $v\in V$ is defined as the number of edges in $E$ which contain $v$.  Two vertices $v$ and $u$ are called \textit{adjacent} (or neighboring) if $(u,v) \in E$. A path between two vertices $u$ and $v$ is an ordered sequence of edges $e_{1}, e_{2}, \ldots, e_{k}$ such that $u \in e_{1}$, $v \in e_{k}$ and for all values of $i$, $e_{i} \cap e_{i+1} \neq \phi$. The graph $G$ is \textit{connected} if there exists a path between any pair of vertices.

To model an Abelian Sandpile Model, we take a connected graph $G$ with a special vertex called \textit{sink}, denoted $s \in V$. Non-sink vertices in $G$ are called ordinary vertices and the set will be denoted by $V_{o} = V - \{s\}$.

\begin{defin}The {\em{configuration}} of a sandpile $G$ is a map $c:V_{o} \rightarrow \mathbb{N}$. It will be represented as a vector.\end{defin}

The configuration $c$ tells us the number of sand particles that each of the ordinary sites currently contain. The \textit{empty} configuration is the zero vector. The \textit{capacity} of a site is the maximum number of particles that it can hold. 

\begin{defin}:An ordinary node $v$ is said to be {\em{unstable}} in a configuration $c$ if $c(v) \geq \text{capacity}(v)$. The configuration $c$ is said to be unstable if any site under it is unstable, else it is referred to as stable. \end{defin}

When a site is unstable it is said to \textit{topple}, i.e. pass on some of its particles to its neighbors. When a site $v$ topples once, it losses $capacity(v)$ particles and each neighbor of $v$ acquires a particle for every edge, common with $v$, that appears in the toppling rule. The sink node never topples. We start with the empty configuration and keep adding particles one by one on sites of our choice and topple when necessary.

The ASM evolves in time through two modes, particle addition at sites and relaxation of unstable sites via toppling. The case of many sites becoming unstable simultaneously also poses no complication as the order in which they are subsequently relaxed does not effect the final stable configuration that is obtained at the end of toppling sequence, hence the prefix abelian. Elementary proofs of such confluence properties can be found in the pioneering paper on ASMs by Dhar \cite{DD90}. A toppling sequence is an ordered set of configurations where every configuration can be obtained from the previous one by toppling some unstable site in it.


\noindent \textbf{Notation:} We write $c_{1} \geq c_{2}$ if $\forall v, c_{1}(v) \geq c_{2}$ and $c_{1} \vdash c_{2}$ if there is a toppling sequence which takes $c_{1}$ to $c_{2}$. Lastly we write, $c_{1} \rightarrow c_{2}$ if $\exists c_{3} \geq c_{1}$ such that $c_{3} \vdash c_{2}$. We say that a configuration $c_{2}$ is \textit{reachable} from $c_{1}$ if $c_{1} \rightarrow c_{2}$ and \textit{unreachable} otherwise. In words, one can add particles to certain sites in $c_{1}$ so that there exists a toppling sequence leading to $c_{2}$. Note that being reachable is a transitive property, i.e. $c_{1} \rightarrow c_{2}, c_{2} \rightarrow c_{3} \Rightarrow c_{1} \rightarrow c_{3}$.

\begin{thm}(\cite{DD90},\cite{auto}) Given any configuration $c$, there exists a unique stable configuration $\sigma(c)$ such that $c \vdash \sigma(c)$, independent of the toppling sequence chosen.
\end{thm}



A configuration is called {\em{recurrent}} if it is reachable from \textit{any} configuration. As already mentioned, we say that a configuration $c_{i}$ is reachable from a configuration $c_{j}$ if by adding some particles to $c_{j}$ and subsequently relaxing it, we can obtain $c_{i}$. A configuration is {\em{transient}} if it is not recurrent. The set of recurrent configurations is therefore, closed under being reachable.

\begin{defin}A configuration $c$ is {\em{recurrent}} iff $\forall c'$ we have $c' \rightarrow c$.\end{defin}

Since recurrence persists under particle addition, we have the following property.

\begin{prop} If configuration $c_{1} \leq c_{2}$, then recurrence of $c_{1}$ implies that of $c_{2}$.\end{prop}

Denote the configuration in which every node $v$ has $\text{capacity}(v) - 1$ particles by $c_{max}$. Clearly, given any stable configuration $c$, one can reach $c_{max}$ simply by adding the required number of particles at each site. Consequently $c_{max}$ is a recurrent configuration. Evidently, every configuration reachable from $c_{max}$ is also recurrent. Checking whether a given configuration is recurrent is an important combinatorial problem. For the sandpile defined above, Dhar \cite{DD90} presents a simple and computationally efficient recursive procedure for deciding recurrence. This algorithm elucidates the importance of the notion of a forbidden sub-configuration. We defer discussing these till next section. Next we introduce the notion of a stochastic abelian sandpile.

A stochastic sandpile (SASM) is the non-deterministic version of the sandpile defined above. Instead of a single toppling rule, a list of toppling rules is associated with each vertex. There is a toppling threshold associated with each node, it is same for each toppling rule, but may not be equal to the vertex degree. When a vertex becomes unstable, a toppling rule is chosen non-deterministically and a toppling, consistent with this chosen rule, occurs.

Formally speaking, an SASM is a triplet $S = (G,T,C)$. Here $G$ is the graph of a deterministic sandpile (with a special sink node) as defined above, $T$ is a function which associates each node with its set of toppling rules, and $C$ is the capacity function which associates with each node its toppling threshold. So $T(v)$ is a collection of multi-sets, such that

\begin{enumerate}
 \item Each member of $T(v)$ has at most $C(v)$ elements.
 \item Each member of $T(v)$ contains elements from the set of ordinary neighbors of $v$, i.e. $V_{o} \cap \text{N}(v)$.
 \item Each edge incident on a vertex $v$ appears in some toppling rule $t \in T(v)$.
\end{enumerate}

When $v$ becomes unstable, one of the members of $T(v)$, say $t$, is chosen and every element in $t$ gets one particle. In case some node appears multiple times in $t$, it gets as many particles. Note that if toppling $t$ takes place at $v$, the number of particles going out of $v$ is $C(v)$ and those received by its neighbors is $|t|$. By convention, the difference between these two numbers is the number of particles that go to sink.

In analogy with the deterministic case, a configuration $c_{f}$ is said to be reachable from a configuration $c_{i}$ if there exists a strategy of adding particles to configuration $c_{i}$ at some sites and choosing toppling rules at each toppling, such that the resulting configuration is $c_{f}$. Consider the special configuration $c_{max}$, in which every site $v$ has $C(v)-1$ particles. Evidently, this configuration is reachable from every other configuration, one just needs to add at every node the number of particles needed to make its weight just less then its toppling threshold. The configuration $c_{max}$ is therefore recurrent. Consequently, all configurations reachable from this $c_{max}$ are recurrent. The ones that are unreachable are transient. This is the definition of recurrence that we will be referring to, throughout this paper.

\begin{defin} A configuration $c$ is recurrent if and only if it is reachable from $c_{max}$.\end{defin}

\section{Forbidden Sub-Configurations and Recurrence}

Consider the following examples of stochastic sandpile. Both of them are based on grid graphs. Each site has stability threshold equal to one, i.e. the sites become unstable if they get two particles. The only difference between them is the toppling rules associated with them.

\begin{itemize}
\item[-] NS-EW (\textit{Manna's Model}): This is the \textit{North South - East West} sandpile. It is an $n\times n$ grid where each site has a threshold of $1$ (that is the site topples when there are two particles on it). Upon toppling it either gives its vertical (north-south) neighbors a particle each or horizontal (east-west) neighbors a particle each. The sites on edges and corners loose some particles to sink depending on the number of edges common with boundary.
\item[-] NE-SW : In analogy with previous example, this is the \textit{North East - South West} sandpile. Threshold of each site is $1$ and when a site topples, either a particle is passed to the northern and eastern neighbors each, or to southern and western ones.
\end{itemize}

\begin{figure}
\centering
\includegraphics[width=.5\textwidth]{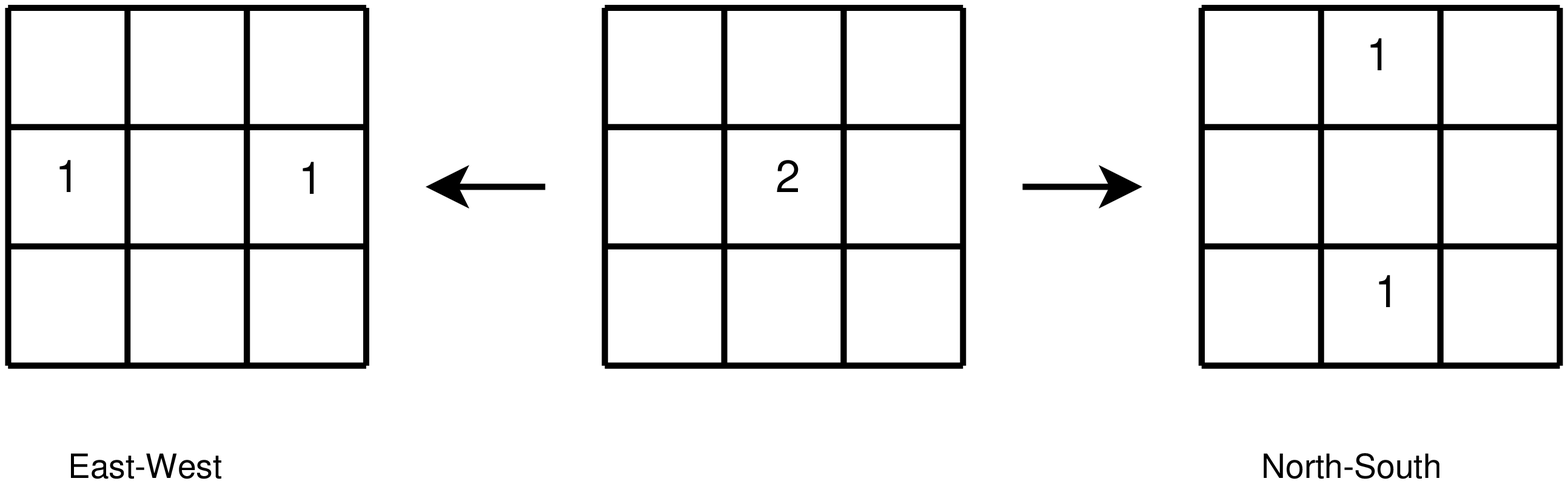}
\caption{NS-EW Stochastic Sandpile}
\end{figure}

\begin{figure}
\centering
\includegraphics[width=.5\textwidth]{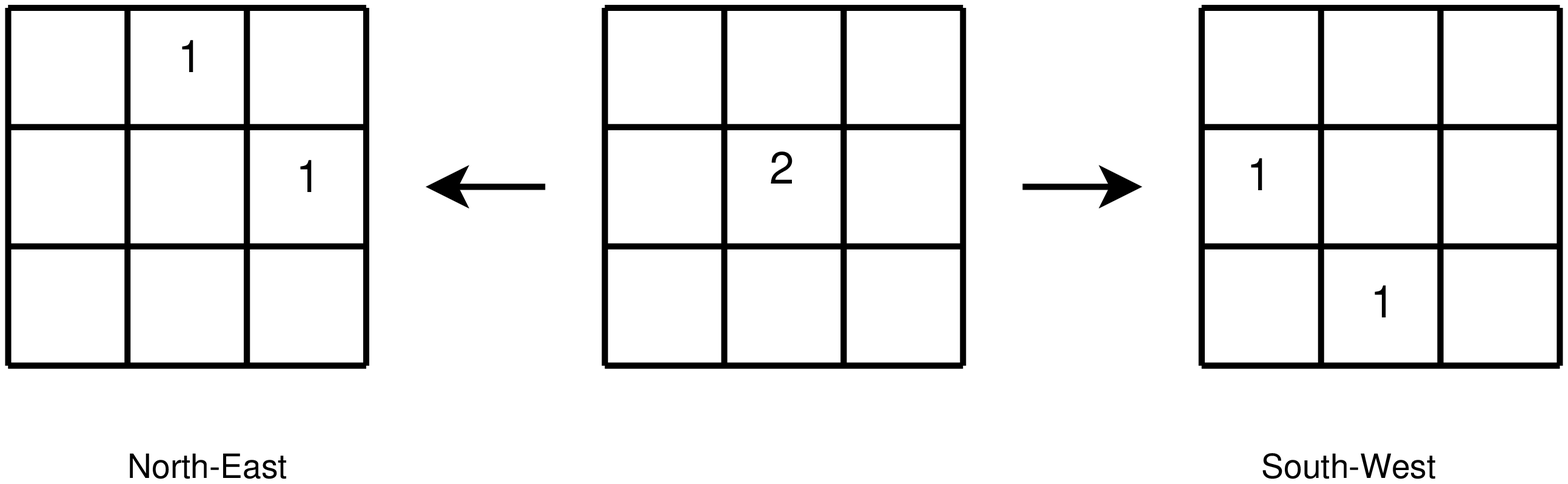}
\caption{NE-SW Stochastic Sandpile}
\end{figure}

The first example, the Manna's model, holds prime importance in being the simplest sandpile which captures the essential features of general stochastic sandpile.

\begin{figure}
\centering
\includegraphics[width=.1\textwidth]{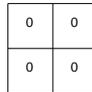}
\caption{A Forbidden Sub-Configuration of the NS-EW sandpile (Manna's model)}
\label{fig:FSC}
\end{figure}

Note that we will often observe only a local region in the sandpile graph. The configuration of that region is a \textit{sub-configuration} of the sandpile under consideration. It so happens that for any given local region, there exists a set of sub-configurations which cannot occur in any recurrent configuration. These are known as forbidden sub-configurations (FSC). We illustrate this with the aid of an example. Consider the Manna's model as defined in the example above. We claim that the sub-configuration shown in figure (\ref{fig:FSC}) consisting of four zeroes is forbidden.

%

The claim that this configuration cannot arise from $c_{max}$ by any sequence of particle additions and toppling is easily verified. For a simple proof, consider the site that toppled last. Depending on which toppling rule was chosen, either the horizontal or the vertical neighbor must have received a particle. It cannot happen that starting from a particle at every node, one manages to reach a state with no particle on any node. More discussion on this property can be found in \cite{DD99}.

Given any configuration, deciding the presence of \textit{any} FSC is equivalent to checking if the given configuration is recurrent. Note that one cannot take existence of forbidden sub-configurations for granted. There exist stochastic sandpile in which all configurations are recurrent and consequently the set of FSCs is empty. The second example defined above, of the NS-EW sandpile, is one such case. It is a simple exercise to demonstrate that, if all sites topple using the same rule then no matter which configuration we start from, the empty configuration is reachable from $c_{max}$. So one infers that there are no FSCs for this sandpile. In general, whether a given stochastic sandpile will have a non-empty set of FSCs is a non-trivial question. In the next section we answer this question. This result is later used in generating arbitrarily large FSCs of Manna's model.

\section{A guarantee for the existence of FSCs}

Consider a sandpile $(G,T,C)$ with $G(V,E)$ as the underlying graph. The following algorithm computes the set of those sites on which some forbidden configuration can exist.

\begin{algo}[\textbf{REDUCE$(G,T,C)$}] The input is a sandpile $(G,T,C)$ and the output is the \textit{reduced} sandpile.
\begin{enumerate}
 \item Initialize $V'$ to $V$ and $F$ to $\phi$.
 \item Select all nodes $v$ in $V'$ such that $\exists$ $t \in T(v)$ for which $t \cap S = \phi$. Call the set of such nodes $D$. If $D$ is empty, then terminate.
 \item Reset $V' = V' - D$ and $F = F \cup D_{s}$.
 \item Repeat step $2$ till $V'$ becomes fixed.
 \item At termination, the set $V'$ contains the sites on which a forbidden configuration can exist and $F$ the set from which all particles can be flushed out. Return the sandpile corresponding to sites in $V$, such that sites in $F$ are deleted from all the toppling rules associated with sites of $V'$ and the threshold function is restricted to sites in $V'$.
\end{enumerate}
\end{algo}

Based on this procedure, we have the following definitions and properties. The algorithm's proof of correctness will follow from these.

\begin{defin}[Irreducible Sandpiles]A sandpile $S = (G,T,C)$ is called \mbox{irreducible} if \textbf{REDUCE($S$)} = \textbf{$S$}
\end{defin}

\begin{defin}[Sub-sandpile] A sandpile $S'$ ($=(G', T', C')$) is called a sub-sandpile of sandpile $S$ ($=(
G, T, C)$) iff

\begin{eqnarray}
V' & \subseteq & V \nonumber \\
C' & = & C \mid_{V'} \nonumber \\
 T'(v) & = &\bigcup_{t \in T(v)} \{ t\cap V' \}  \quad \forall v \in V' \nonumber
\end{eqnarray}

we will also denote this sandpile by $S\mid_{V'}$.
\end{defin}

\begin{prop}
 A sub-sandpile is a valid sandpile, i.e. particles are not created during toppling and sink is reachable from each site.
\end{prop}

A sub-sandpile is a restriction of the original sandpile to some subset of sites. A special case is the sandpile obtained by \textit{deleting} a site $v$. Here, $v$ is removed from the set of sites and also from all the toppling rules it appears in. Denote this sandpile by $S\mid_{V-v}$. The idea behind the notion of irreducibility is that while checking for recurrence of a certain configuration, certain sites may not contribute to the problem complexity at all. Overlooking these can reduce the problem size substantially. The sandpiles which are, structurally speaking, the simplest to analyze, are defined below.

\begin{defin}[Minimal Irreducible Sandpiles] A sandpile $S$ $=(G,T,C)$ is called \textit{minimal irreducible} if and only if $\forall v \in V$, \textbf{REDUCE($S\mid_{V-v}$)} = \textbf{$\phi$}
\end{defin}

That is deleting \textit{any} site from the sandpile makes its analysis trivial. Following are some useful properties of irreducible sandpiles.

\begin{prop}
For a sandpile $S$ containing a site $v$, there may not be any minimal irreducible sub-sandpile containing $v$.
\end{prop}

\begin{figure}
\centering
\includegraphics[width=.5\textwidth]{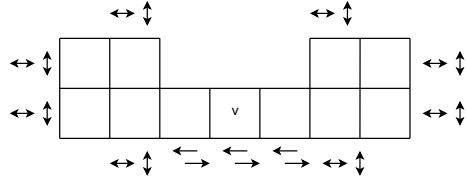}
\caption{An example containing a site which does not belong to any minimal irreducible sub-sandpile}
\label{fig:irredASM}
\end{figure}

As an example, consider the sandpile in figure \ref{fig:irredASM}. Here, all the sites are shown with the list of toppling rules associated with them. The sites with single headed arrows have unit threshold, and those with double sided arrows have threshold $2$. It is clear that no matter which site we remove, the \textbf{REDUCE} procedure will leave only one or both the square blocks at the end, which means there doesn't exist any minimal irreducible sub-sandpile which contains $v$.

\begin{prop}\label{prop:emptyFSC}
 The empty configuration on any irreducible sandpile is forbidden.
\end{prop}

\textbf{Proof of Property \ref{prop:emptyFSC} : }The proof of this is by contradiction, it is a simple extension of the idea behind demonstrating the smallest FSC of Manna's model in previous section. Consider the last site which toppled. No matter which toppling rule was used, there must be a site in the sandpile which should have received particle as a consequence (since the sandpile is irreducible). Hence the configuration containing all zeroes cannot be achieved via toppling and particle additions to $c_{max}$. $\square$

The structure of the sandpile (in context of reducibility) and its weight (total particle count) behave in opposite directions with regards to the recurrence of a configuration. Property \ref{prop:emptyFSC} poses a lower limit on the structural complexity of a sandpile when the weight is zero, to ensure transience. The next property is complementary to the previous one, in that it poses an upper limit on the weight for a given structural complexity.

\begin{prop}\label{prop:minmax}
The configuration with one particle on any site of a minimal irreducible sandpile is recurrent.
\end{prop}

To prove this, we will show that this configuration is obtainable from a configuration in which all sites are unstable and is consequently reachable from $c_{max}$. Let $S = (G,T,C)$ be the irreducible sandpile and $v$ be any given vertex. Run the algorithm REDUCE on the sub-sandpile $S|_{V-{v}}$. Since $S$ is minimal irreducible, the algorithm will terminate with a null set. Denote the set removed from $S$ in the $i^{th}$ iteration by $D_{i}$. Then $D_{i}$s have the following property.

\begin{prop}
 $\forall v_{j} \in D_{i}, \exists t \in T(v_{j})$ such that $t \cap D_{k} = \phi $ for all $k \geq i$.
\end{prop}

This means that each of the vertices removed in $i^{th}$ iteration have at least one toppling rule which affects only those vertices which have been removed in preceding iterations. The validity of this property follows from the definition of algorithm trivially.

\textbf{Proof of Property \ref{prop:minmax} : } For every vertex $v_{i}$ define the sets $D_{i,1}$, $D_{i,2}$, \ldots , $D_{i,k_{i}}$, where $k_{i}$ is the number of iterations taken by the algorithm to reduce the sub-sandpile to null, as above. Now consider the configuration $c$ in which only one site, say $v_{m}$, has one particle. All other sites are vacant.

\begin{eqnarray}
 c(v_{j}) = \delta_{m,j} \nonumber
\end{eqnarray}

here, $c(v_{j})$ denotes the number of particles at $v_{j}$. $\delta_{i,j}$ is the \textit{Kroneckor delta function}, its value is $1$ when $i=j$ and zero in all other cases. We will now show that,

\begin{enumerate}
 \item $c$ is obtainable from a configuration $c$' in which all sites in $D_{m,1}$ can have arbitrarily large heights
 \item Any $c$' in which sites in $D_{m,1}$, $D_{m,2}$, \ldots, $D_{m,k}$ have arbitrarily large heights is obtainable from a configuration $c$'' in which all sites in $D_{m,1}$, $D_{m,2}$, \ldots, $D_{m,k}$, $D_{m,k+1}$ have arbitrarily large heights
\end{enumerate}

 These two properties would imply that $c$ is obtainable from a configuration in which all sites have arbitrarily large heights. We first prove the first part. $c_{1}$ has one particle at $v_{m}$ and zero elsewhere. Without loss of generality assume that $m=1$.
\begin{eqnarray}
 c_{1}(v_{j}) = \delta_{1,j} \nonumber
 \end{eqnarray}

By minimal irreducibility of $S$, $\mid D_{1,1} \mid \geq 1$. Choose any member of $D_{1,1}$, say $v_{2}$. The definition of $D_{1,1}$ implies $\exists t \in T(v_{2})$ such that $ v_{1} \in t $ . Therefore, using this toppling rule, $c_{1}$ is reachable from the following $c_{2}$,

\begin{eqnarray}
 c_{2}(v_{j}) = C(v_{2}).\delta_{2,j} \nonumber
 \end{eqnarray}

On similar lines, let $v_{3}$ be an element of $D_{2,1}$, then $c_{2}$ is reachable from the following $c_{3}$,

\begin{eqnarray}
 c_{3}(v_{j}) = (C(v_{2})-1).\delta_{2,j} + C(v_{3}).\delta_{3,j} \nonumber
\end{eqnarray}

In general, $c_{k}$ is reachable from the following $c_{k+1}$,

\begin{eqnarray}
 c_{k+1}(v_{j}) = \sum_{i<k}c(v_{i}).\delta_{i,j}  + (c(v_{k})-1).\delta_{k,j} + C(v_{k+1}).\delta_{k+1,j} \nonumber
\end{eqnarray}

here $v_{k+1} \in D_{k,1}$. The finiteness of $G$ implies that some $v_{k}$ is revisited at some step $m$. So,

\begin{eqnarray}
 c_{m}(v_{j}) = \sum_{i<m}(C(v_{i}) - 1).\delta_{i,j}  + C(v_{k}).\delta_{k+1,j} \nonumber
\end{eqnarray}

Here, $v_{m}$ is the same as $v_{k}$ that is the loop is formed at vertex $v_{m}$. From the corresponding configuration $c_{m}$, repeat the whole procedure again, except that now at each step, choose from those members of the $D_{l,1}$s under consideration, which have been chosen the least number of times. In $max_{l}\mid D_{l,1} \mid$ iterations, every vertex in every $D_{l,1}$ will have at least $C(v) - 1$ particles. By repeating the whole process $3.max_{l}\mid D_{l,1} \mid$ times, every node in every $D_{l,1}$ will become unstable.
We have assumed that for all nodes $L(v) > 1$, which is a realistic assumption, since every site in the sandpile has non-zero capacity.

 \begin{eqnarray}
 L(v) \geq 2 \Rightarrow 3.(L(v) - 1) \geq L(v) \nonumber
 \end{eqnarray}

This completes the proof of the first part. The proof for second part follows the same idea. Let $c$ be a configuration in which nodes from $D_{m,1}$, $D_{m,2}$, \ldots, $D_{m,l}$ have arbitrarily large heights. Then $c$ is obtainable from the following $c_{1}$,

\begin{center}
$\forall v_{i} \in D_{m,l+1} \exists t \in T(v_{i})$ such that\\
$t \cap D_{m,l+j} = \phi \quad  \forall j \geq 1 $ \\
\end{center}

choose one such $t$ for every site $v_{i}$ in $D_{m,l+1}$ and denote it $t_{i}$. Let the size of $D_{m,l+1}$ be $r$. Then the configuration is obtainable from a configuration $C'$ in which all sites in $D_{m,l+1}$ are unstable and heights of sites in $D_{m,1}$, $D_{m,2}$, \ldots, $D_{m,l}$ is reduced by at most $r.(max_{v_{i} \in D_{m,l+1}} L(v_{i}))$. Thus the second part is also proved. So, \textit{for every minimal irreducible sandpile, any configuration with just one particle on any site is recurrent}. $\square$

\begin{defin}[Weight Maximal Transient Configurations]
A configuration $c$ over a sandpile $S$ is weight maximal transient if it is transient and adding any single particle to it at any site makes it recurrent.
\end{defin}

\begin{cor}
The empty configuration over a minimal irreducible sandpile is weight maximal transient.
\end{cor}

\noindent \textit{Remark:} The correctness of procedure \textbf{REDUCE} can be proved using essentially the same principles, which imply the properties elucidated above. The procedure terminates with a sub-sandpile such that for every site in it, every toppling rule effects some other site in the sub-sandpile. The empty configuration on such a sandpile is transient. This can be shown if one considers the last toppling before the empty configuration was attained. No matter which site toppled using whichever toppling rule, some other site must have received a particle. Hence, empty configuration is not reachable. The procedure therefore produces a set of sites, on which the empty configuration is transient.

\section{Generating FSCs for Manna's model}

In this section we will outline a scheme to generate arbitrarily large FSCs for the Manna's model (the NS-EW stochastic sandpile defined previously). Dhar \cite{DD99} discusses the FSCs of Manna's model in detail and mentions some of them. These FSCs were derived empirically by running extensive simulations. We will outline a scheme which efficiently generates most of these.

\begin{prop}\label{prop:comp1}
If $H_{1}$ and $H_{2}$ are irreducible sub-sandpiles, then
\begin{enumerate}
 \item $H_{1} \cup H_{2}$ is irreducible
 \item Any configuration over $H_{1} \cup H_{2}$ with just one particle is forbidden.
 \item (Gluing Property of forbidden configurations) If the following conditions hold
 \begin{enumerate}
 \item $c_{1}$ is a forbidden configuration over $H_{1}$ and $c_{2}$ forbidden over $H_{2}$
 \item $c_{1}$ and $c_{2}$ agree over the intersection of $H_{1}$ and $H_{2}$ that is $c_{1}(v) = c_{2}(v) \forall v \in V(H_{1}) \cap V(H_{2})$
 \end{enumerate}
 denote by $c_{1} \vee c_{2}$ the configuration on $V(H_{1}) \cup V(H_{2})$ which agrees with $c_{1}$ on $V(H_{1})$ and $c_{2}$ on $V(H_{2})$. Then, the configuration obtained by adding an extra particle to $c_{1} \vee c_{2}$ at any $v_{k} \in V(H_{1}) \cap V(H_{2})$  such that $\nexists t \in T(v_{k})$ with $t = \{ v_{i} \}$ is also forbidden. Basically, we are adding a particle at a site in $V(H_{1}) \cap V(H_{2})$ which is not connected to any other site in $V(H_{1}) \cap V(H_{2})$.
\end{enumerate}
\end{prop}

\textbf{Proof of Property \ref{prop:comp1} : } The proves are as follows
\begin{enumerate}
 \item Since both $H_{1}$ and $H_{2}$ are irreducible, any site in $H_{1} \cup H_{2}$ belongs to at least one of the component sandpiles and hence has at least one member from every toppling rule in the composite sandpile. Therefore, it is irreducible.
 \item Empty configurations over each of the sub-sandpile is forbidden (property \ref{prop:emptyFSC}). The given configuration over the union contains at least one of these as a sub-configuration. Hence it is also forbidden.
 \item In all configurations which can generate a configuration satisfying property \ref{prop:comp1} by one toppling, the particle must come from some site in either $H_{1}$ or $H_{2}$. All these configurations, therefore must contain either $c_{1}$ or $c_{2}$ as a sub-configuration (over $H_{1}$ or $H_{2}$, as the case may be); both of which are forbidden. Hence, the configuration under consideration is also forbidden.
\end{enumerate} $\square$

\begin{figure}
\centering
\includegraphics[width=.7\textwidth]{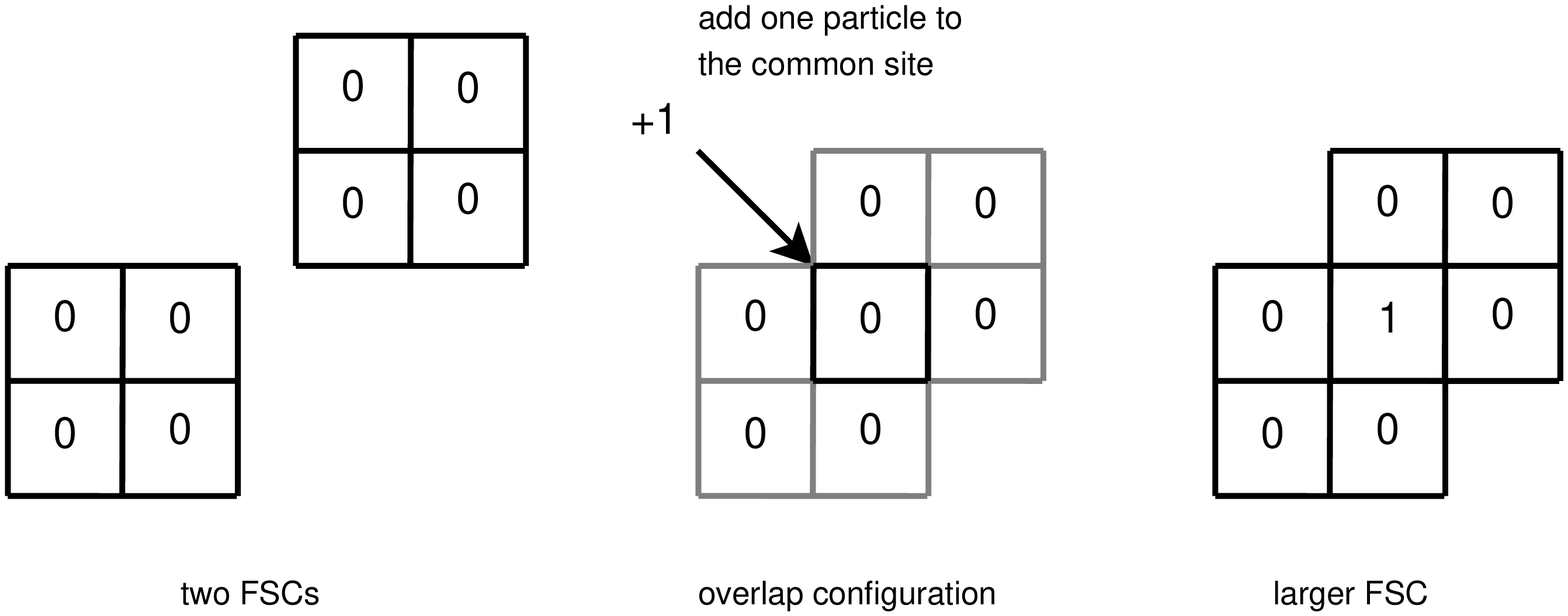}
\caption{A Forbidden Sub-Configuration of the Manna's model}
\label{fig:FSCover-1}
\end{figure}

\begin{figure}
\centering
\includegraphics[width=.7\textwidth]{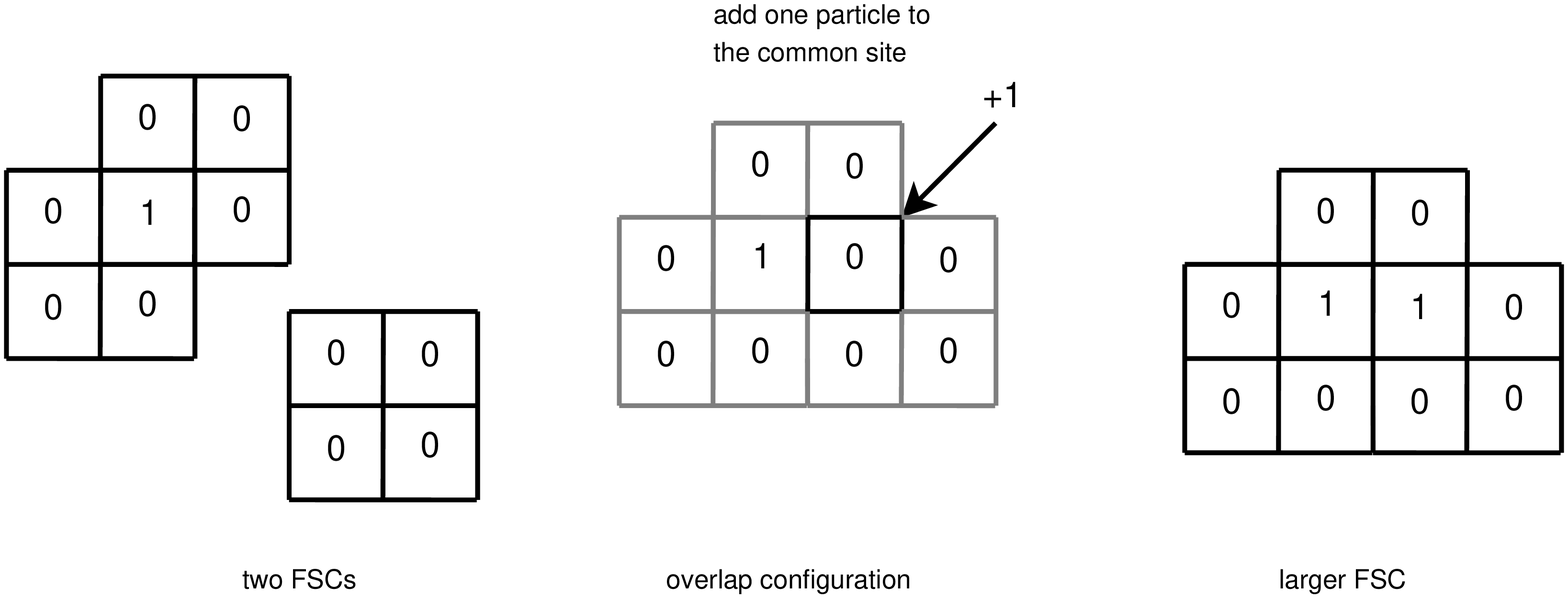}
\caption{A Forbidden Sub-Configuration of the Manna's model}
\label{fig:FSCover-2}
\end{figure}

Using this property, we will now proceed with outlining our procedure for generating FSCs for the Manna's model. Recall that this is just the grid based sandpile where each site has a unit threshold and can topple \textit{vertically} (the vertical neighbors get a particle each) or \textit{horizontally} (the horizontal neighbors get a particle) each. The smallest FSC is shown in figure (\ref{fig:FSC}). We will use two copies of this to create a larger FSCs using the property \ref{prop:comp1}. See figure \ref{fig:FSCover-1}. The two smaller configurations are arranged so that they have a common \textit{corner} site. An additional particle at this site does not make the configuration recurrent. It is easy to check that these configurations and there union satisfy all the preconditions required by property \ref{prop:comp1}. This gives us our first FSC, which is not minimal. Repeated application of this rule, gives us a yet larger FSC. See figure \ref{fig:FSCover-2}.

\section{Future Work}

The main open question is that of coming up with a polynomial time recurrence checking algorithm. A randomized procedure with some performance guarantee will also be a highly significant. However we believe that the logical next step from here is coming up with an efficient scheme to generate all the forbidden configurations. There is exactly one FSC in Dhar's list \cite{DD99} which our procedure can not generate. It is not known if there are other configurations like this. Improving the scheme we have proposed to account for these missing configurations is going to be a crucial step in tackling the original problem.

\noindent \textbf{Acknowledgements} We thank Prof. Deepak Dhar (TIFR, Mumbai) for introducing us to this problem. We also thank Prof. Sundar Vishwanathan (IIT Bombay, Mumbai) for the guidance throughout the course of this work.

\end{document}